\documentclass[11pt]{article}

\usepackage[T1]{fontenc}
\usepackage[latin1]{inputenc}

\usepackage{amssymb}
\usepackage{amsmath}
\usepackage{stmaryrd}
\usepackage{array}
\usepackage{hhline}
\usepackage{ulem}
\usepackage{epsfig}
\usepackage{psfrag}
\usepackage{afterpage}
\usepackage{theorem}

\usepackage{eufrak}

\newlength{\listeespacee}
\newlength{\tableautasse}
\newcommand{\dd}{\text{\rm d}}
\newcommand{\SSS}{{\cal S}}

\newcommand{\NN}{{\cal N}}
\newcommand{\UU}{{\cal U}}

\newcommand{\VV}{{\cal V}}

\newcommand{\II}{{\cal I}}

\newcommand{\KK}{{\cal K}}

\newcommand{\FF}{{\cal F}}

\newcommand{\zbar}{{\overline{z}}}

\newcommand{\R}{\mathbb R}
\newcommand{\N}{\mathbb N}
\newcommand{\Z}{\mathbb Z}
\newcommand{\C}{\mathbb C}
\newcommand{\K}{\mathbb K}

\newcommand{\Id}{\operatorname{Id}}

\newcommand{\End}{\operatorname{End}}

\newcommand{\Mat}{\mbox{\rm Mat}}

\newcommand{\im}{\operatorname{Im}}

\renewcommand{\Im}{\operatorname{Im}}

\theoremstyle{change}
\newtheorem{enonce}{}[section]

\newtheorem{de}[enonce]{Definition}
\newtheorem{prop}[enonce]{Proposition}

\newtheorem{add-property}[enonce]{Additional property}
\newtheorem{lem}[enonce]{Lemma}

\newtheorem{te}[enonce]{Theorem}
\newtheorem{cor}[enonce]{Corollary}

\theorembodyfont{\rmfamily}
\newtheorem{rem-notation}[enonce]{Remark/Notation}
\newtheorem{def-notation}[enonce]{Definition/Notation}
\newtheorem{notation}[enonce]{Notation}
\newtheorem{rem}[enonce]{Remark}
\newtheorem{example-rem}[enonce]{Example/Remark}

\newtheorem{consequence-notation}[enonce]{Consequence/Notation}

\newtheorem{recall-rem}[enonce]{Recall/Remark}

\newtheorem{recall/vocabulary}[enonce]{Recall/vocabulary}
\newtheorem{vocabulary}[enonce]{Vocabulary}

\newtheorem{example}[enonce]{Example}

\newcounter{claimcounter}
\newlength{\aux}

\psfrag{S}{${\cal S}=\sigma(\overline{{\cal V}})$}
\psfrag{Vbar}{\fbox{$\overline{{\cal V}}=\pi({\cal V})$}}
\psfrag{V}{\fbox{${\cal V}$}}
\psfrag{coordinates zbar}{coordinates $\zbar_{i,j}$ on $\overline{\cal V}$}
\psfrag{vectors Z}{vectors $Z_{i,j}=\dd\sigma.\frac{\partial}{\partial \zbar_{i,j}}$ (in black, along ${\cal S}$)}
\psfrag{vectors Z0}{vectors $Z^{(0)}_{i,j}$ (in grey, everywhere in ${\cal V}$)}
\psfrag{pi}{$\pi$}
\psfrag{leaves of I}{leaves of $\cal I$}
\psfrag{sigma}{$\sigma$}

\begin{document}
\normalem
\begin{center}
{\Large\bf An integrability condition for fields\\of nilpotent endomorphisms\bigskip\\

}Charles Boubel\footnote{Institut de Recherche Mathématique Avancée, UMR 7501 -- 
Université de Strasbourg et CNRS, 7 rue René Descartes, 67084 STRASBOURG CEDEX, FRANCE}, March 28th. 2011.
\end{center}

\setlength{\aux}{\textwidth}\addtolength{\aux}{-2cm}\noindent\hspace*{\fill}{\small\begin{tabular}{p{\aux}}{\bf Abstract.} We give a necessary and sufficient condition on the $1$-jet of a field of nilpotent endomorphisms to be integrable. Together with the well known corresponding condition for an almost complex structure, the nullity of its Nijenhuis tensor, this gives an integrability condition for any field of endomorphisms.\medskip\end{tabular}}\hspace*{\fill}\medskip\\

It is a classical question to ask whether (the germ of) an almost complex structure $J$ is (the germ of) a complex structure {\em i.e.\@} if it is integrable: does it exist local coordinates in which $J$ becomes a constant matrix, namely {\tiny$\left(\begin{array}{cc}0&I\\-I&0\end{array}\right)$}? A well known necessary and sufficient condition on the $1$-jet of $J$ is that the Nijenhuis torsion tensor ${\cal N}_J$ of $J$ vanishes \cite{newlander-nirenberg}. We address here the same question for a field of nilpotent endomorphisms $A$: instead of ``$J^2=-\Id$'', we take ``$A^n=0$'' for some $n$. More precisely, we suppose that $A$ is conjugate, at every point, to some fixed nilpotent endomorphism ---~this is necessary to hope integrability.

Immediately, the nullity of the Nijenhuis tensor ${\cal N}_A$ is necessary. Simple examples show that this condition is not sufficient, see section \ref{bonus}, see also \cite{lavandier2}. We show here that it becomes sufficient, together with the additional condition that each distribution of the flag $(\ker A^p)_{p=1}^{n-1}$ is involutive. The proof, unlike that of \cite{newlander-nirenberg}, follows essentially from the Cauchy-Lipschitz theorem and some standard differential calculus.

In combination with the integrability condition for complex structures, this immediately gives an integrability condition for any smooth endomorphism field $M$: $M$ is integrable if and only if it has constant invariant factors, ${\cal N}_M=0$ and $\ker(P(M))$, for each invariant factor $P$ of $M$, is involutive. 

A general viewpoint on this type of problems, that we do not use here, is given in \cite{gromov_rtg}.\medskip

\noindent{\bf Thanks.} I thank R.\@ Bryant, T.\@ Delzant, \'E.\@ Ghys, A.\@ Oancea, and E.\@ Opshtein for their quickly answering my questions.\medskip

Everywhere, $A$ is a germ of endomorphism field of $T\K^d$ around $0$ in $\K^d$, with $\K=\R$ or $\K=\C$, {\em i.e.\@} a smooth (holomorphic if $\K=\C$) section of $\End_\K(T\K^d)$ on a neighbourhood $\cal V$ of $0$. All objects: coordinates, tangent bundles {\em etc.} are real if $\K=\R$ and complex if $\K=\C$.

Section \ref{resultats} recalls the requisite definitions and states the results, section \ref{preuves} gives the proofs and section \ref{bonus} provides some additional results, comments and examples.

\section{Definitions and results}\label{resultats}

We recall the two following definitions.

\begin{de}\label{def_torsion}
The Nijenhuis torsion tensor of $A$ is the vector valued $2$-form defined by:
$${\cal N}_A(X,Y)=[AX,AY]-A[X,AY]-A[AX,Y]+A^2[X,Y].$$
We let the reader check that it is a tensor, see {\em e.g.\@} \cite{kono1}, ch.\@ 1 prop.\@ 3.12, where the torsion tensor $S_{A,B}$ of some couple $(A,B)$ of fields of endomorphisms is introduced. Our ${\cal N}_A$ is equal to $\frac12S_{A,A}$. The fact that $\K=\R$ or $\K=\C$ plays no role here.
\end{de}

\begin{de}
The field $A$ is called integrable if there exists, on a neighbourdood $\cal V$ of the origin, a coordinate system in which $\Mat(A)$ is constant {\em i.e.\@} a diffeomorphism, or a biholomorphism $\varphi:{\cal V}\rightarrow {\cal U}\subset\K^d$ such that $\varphi_\ast A$ is the restriction to $\cal U$ of a linear transformation of $\K^d$.
\end{de}

Here we show the following result.

\begin{te}\label{theorem} Let $A$ be a germ of field of nilpotent endomorphisms of order $n\geqslant 1$ on $\K^d$. If $\K=\C$, we take $A$ holomorphic. If $\K=\R$, we take $A$ of class $C^{\omega}$, $C^{\infty}$ or $C^{r}$ with $r\geqslant n-1$.

Then $A$ is integrable if and only if the three following conditions are satisfied:\medskip

-- the invariant factors of $A$ are constant,\medskip

-- ${\cal N}_A=0$,\medskip

-- each distribution $\ker A^p$, for $p\in\N$, is involutive ---~hence integrable.\medskip

\noindent If $A$ is analytic or of class $C^\infty$, the integral coordinates have the same regularity. If $A$ if of class $C^r$ with $n-1\leqslant r<\infty$, they are at least, and possibly not more than, of class $C^{r-n+2}$. If $A$ satisfies the three conditions but is not of class $C^{n-1}$, it is non integrable in general.
\end{te} 

\begin{rem}The regularity condition ``class $C^{n-1}$'', though minor, has to be mentioned. In other equivalence problems of $G$-structures of order 1 (see \cite{gromov_rtg}), with $G$ reductive, and solved as P.D.E.\@ problems, such a strong regularity condition does not seem to appear (see {\em e.g.\@} \cite{libermann} or Theorem II of \cite{nijenhuis-woolf}). In th.\@ \ref{theorem}, the group $G$ is the centraliser of ${\rm Id} + A$ in GL$_d(\K)$, which is not reductive as soon as $A\neq 0$. The regularity condition seems to be linked to that fact.

The present coordinates are not the solution of an elliptic P.D.E., see Remark \ref{non_elliptique}. Instead, they arise naturally as the solution of O.D.E.'s integrated by induction. In that sense the proof of Theorem \ref{theorem} is similar to that of the Frobenius criterion given in \cite{hormander}, C.1.1.
\end{rem}

Together with the classical integrability condition for complex structures, the present result gives easily the following corollary.

\begin{cor}\label{corollaire}
If $A$ is any field of endomorphisms of class $C^{\infty}$ on $\R^d$, with constant invariant factors, then $A$ is integrable if and only if the three following conditions are realised:\medskip

-- the invariant factors of $A$ are constant,\medskip

-- $\NN_A=0$,\medskip

-- the distribution $\ker(P(A))$, for each invariant factor $P$ of $A$, is involutive.
\end{cor}

\noindent Of course, the minimal regularity condition in general is that $A$ is of class $C^{n-1}$ along each integral leaf of $\ker(P(A))$, with $P=Q^n$, $Q$ irreducible. Eventually, a little remark, proven in section \ref{preuves}, is worth to be pointed out autonomously.

\begin{rem}\label{images} If $A$ is nilpotent, the nullity of ${\cal N}_A$ implies the integrability of each distribution $\Im A^p$, but not that of the kernel distributions $\ker A^p$.
\end{rem}

\section{Proof of the results}\label{preuves}

If $A$ is integrable, it is conjugate, at any point, to some fixed nilpotent matrix, so it has constant invariant factors. So the first condition of Theorem \ref{theorem} and of Corollary \ref{corollaire} is the $0$-order integrability condition for $A$, and is necessary. From now on we suppose it holds.\medskip

We introduce the following technical torsion-related tensor, and one of its properties.

\begin{de}
If $B$ is another endomorphism field on $\VV$ and if $A$ and $B$c ommute, we introduce:
$${\cal N}'_{A,B}(X,Y)=[AX,BY]-A[X,BY]-B[AX,Y]+AB[X,Y].$$
The reader may check it is a tensor; the sum $S_{A,B}={\cal N}'_{A,B}+{\cal N}'_{B,A}$ is the torsion of $(A,B)$ cited in Def.\@ \ref{def_torsion}, well-defined even if $AB\neq BA$. So here ${\cal N}_A={\cal N}'_{A,A}=\frac12S_{A,A}$.
\end{de}

\begin{prop}\label{technique}
All ${\cal N}'_{A^p,A^q}$ for $p,q\in\N^\ast$ depend only on ${\cal N}_{A}$, through both following relations:\medskip

{\bf (i)} for all couple $(X,Y)$ of vectors, ${\cal N}'_{A,A^q}(X,Y)=\sum_{k=1}^{q}A^{q-k}{\cal N}_A(X,A^{k-1}Y)$,\medskip

{\bf (ii)} for all couple $(X,Y)$ of vectors, ${\cal N}'_{A^p,A^q}(X,Y)=\sum_{k=1}^{p}A^{p-k}{\cal N}'_{A,A^q}(X,A^{k-1}Y)$.\medskip

\noindent In particular, if ${\cal N}_A=0$, then all ${\cal N}'_{A^p,A^q}$ and all ${\cal N}_{A^p}={\cal N}'_{A^p,A^p}$ also vanish.
\end{prop}

\noindent{\bf Proof.} As $\NN'_{A,A^q}(X,Y)=-\NN'_{A^q,A}(Y,X)$, {\bf (i)} is a special case of {\bf (ii)}. Let us prove {\bf (ii)} by induction on $p$. It is trivial for $p=1$. Suppose it holds for some $p$.
\begin{align*}
{\cal N}'_{A^{p+1},A^q}(X,Y)&=[A^{p+1}X,A^qY]-A^q[A^{p+1}X,Y]-A^{p+1}[X,A^qY]+A^{p+q+1}[X,Y]\\
&=[A^{p+1}X,A^qY]-A^{q}[A^{p+1}X,Y]-A[A^{p}X,A^qY]+A^{q+1}[A^{p}X,AY]\\
&\phantom{=\quad}+A[A^{p}X,A^qY]-A^{q+1}[A^{p}X,AY]-A^{p+1}[X,A^qY]+A^{p+q+1}[X,Y]\\
&={\cal N}'_{A,A^q}(A^{p}X,Y)+A\left({\cal N}'_{A^p,A^q}(X,Y)\right),
\end{align*}
hence it holds for $p+1$.\hfill$\Box$\bigskip

\noindent{\bf Proof of Remark \ref{images}.} Now we can prove Remark \ref{images}. As ${\cal N}_A=0$, each distribution $\Im A^p$ is integrable. Let us prove it is involutive, the conclusion follows by the Frobenius criterion. Let us take $X$ and $Y$ any vector fields and show: $[A^pX,A^pY]\in \Im A^p$. By Proposition \ref{technique}, ${\cal N}_{A^p}(X,Y)=0$, so $[A^pX,A^pY]=A^p[X,A^pY]+A^p[A^pX,Y]-A^{2p}[X,Y]$ and we are done. Besides, example \ref{noyau_non_integrable} gives a counter example to the integrability of $\ker A^p$.

\begin{notation} If $A$ satisfies the three conditions of Theorem \ref{theorem}, using Remark \ref{images}, we denote respectively by ${\cal I}^p$ and ${\cal K}^p$ the integral foliation of the distribution $\Im N^p$, respectively $\ker N^p$, for any $p$. 
We shortly denote $\II^1$ by $\II$, and denote by $\pi$ the projection $\VV\rightarrow\VV/\II$.
\end{notation}

\begin{recall/vocabulary} If $\FF$ is some foliation on $\VV$, and $V$ some distribution or vector field on $\VV$, $V$ is called {\em basic} (for $\FF$) if, for any open set $\UU$ where $\FF$ is trivial, setting $\Pi:\UU\rightarrow\UU/\FF$, $\Pi_\ast V$ is constant along each leaf of $\FF$, and so $V$ ``passes to the quotient'' on $\UU/\FF$. If $V$ is a vector field, this means exactly that its flow sends each leaf of $\FF$ on a leaf of $\FF$.
\end{recall/vocabulary}

To prove the theorem, we already introduce the following, and prove a lemma about it.\medskip

The two flags $\ker A\subset\ker A^2\subset\ldots\subset\ker A^{n-1}\subset\ker A^n=T\K^d$ and $\Im A^{n-1}\subset\ldots\subset\Im A\subset\Im A^0=T\K^d$ satisfy the following inclusion properties:
$$\begin{array}{cccccccccc}
\Im A^{n-1}\\
\cap\\
(\Im A^{n-2}\cap\ker A)&\subset&\Im A^{n-2}\\
\cap&&\cap\\
(\Im A^{n-3}\cap\ker A)&\subset&(\Im A^{n-3}\cap\ker A^2)&\subset&\Im A^{n-2}\\
\cap&&\cap&&\cap&&\ddots\\
\vdots&&\vdots&&\vdots\\
\ker A&\subset&\ker A^2&\subset&\ker A^3&\subset&\cdots&\ker A^{n-1}&\subset&T\VV.\\
\end{array}$$
Any entry of this array is integrable, generating the following foliations:
$$\label{tableau_feuilletages}\begin{array}{cccccccccc}
\II^{n-1}\\
\cap\\
(\II^{n-2}\cap\KK^1)&\subset&\II^{n-2}\\
\cap&&\cap\\
(\II^{n-3}\cap\KK^1)&\subset&(\II^{n-3}\cap\KK^2)&\subset&\II^{n-2}\\
\cap&&\cap&&\cap\\
\vdots&&\vdots&&\vdots\\
\KK&\subset&\KK^2&\subset&\KK^3&\subset&\cdots&\KK^{n-1}&\subset&\VV,\\
\end{array}$$

\begin{lem}\label{section} If the three conditions of Theorem \ref{theorem} are realised, then there exist local coordinates $((x^{i,j}_\alpha)_{\alpha})_{n\geqslant i\geqslant j\geqslant1}$ adapted to this array of foliations {\em i.e.\@} such that, for any $p\in\llbracket1,n\rrbracket$ and $q\leqslant p$, the $(x^{i,j}_\alpha)_{\alpha}$ with $i\leqslant p$ and $j\leqslant q$ parametrise the leaves of $\II^{n-p}\cap\KK^{q}$. 
The coordinates may be chosen of class $C^{r+1}$ in case $\K=\R$, holomorphic in case $\K=\C$.
\end{lem}


\noindent{\bf Proof.} The lemma is nothing but the fact {\bf\mathversion{bold}$(\ast)$\mathversion{normal}} that the distributions $\ker A^q$ are basic for any of the foliations $\II^p$, or in other words, that the foliations $\KK^q$ ``pass'' to the quotient by any of the $\II^p$. Indeed if {\bf\mathversion{bold}$(\ast)$\mathversion{normal}} holds, take any coordinate system $(x^i_\alpha)_{i=1}^n$ such that the leaves of $\II^{n-p}$ are the levels of $((x^i_\alpha)_\alpha)_{i>p}$. In particular, $(x^n_\alpha)_\alpha$ may be viewed as coordinates of $\pi(\VV)$. By {\bf\mathversion{bold}$(\ast)$\mathversion{normal}}, $\pi(\VV)$ is endowed with the foliations $\pi(\KK)\subset\pi(\KK^2)\subset\ldots\subset\pi(\KK^{n-1})$, so $(x^n_\alpha)_\alpha$ may be turned into some other system $((x^{n,1}_\alpha)_\alpha,\ldots,(x^{n,n}_\alpha)_\alpha)$, adapted to this flag: the leaves of $\pi(\KK^q)$ are the levels of $((x^{n,i}_\alpha)_\alpha)_{i>q}$. Inductively, we build the coordinates of the lemma.

Now, {\bf\mathversion{bold}$(\ast)$\mathversion{normal}} amounts exactly to a stronger version of Remark \ref{images}: any of the distributions $\ker A^p+\Im A^q$ is involutive. We prove it and are done. Take $X,X'$ vector fields in $\ker A^p$ and $Y,Y'$ two vector fields in $\Im A^q$. Then: $[X+Y,X'+Y']=[X,X']+[Y,Y']+[X,Y']+[Y,X']$. As $\ker A^p$, by assumption, and $\Im A^q$, by Remark \ref{images}, are involutive, $[X,X']\in\ker A^p$ and $[Y,Y']\in\Im A^q$. We are left with showing, for instance, that $[X,Y']\in\ker A^p+\Im A^q$ {\em i.e.\@} that $A^p[X,Y']\in\Im A^{p+q}$. Take a field $Z$ such that $Y'=A^qZ$:
$$A^p[X,Y']=-\underbrace{{\cal N}'_{A^p,A^q}(X,Y)}_{=0\text{ by Prop \ref{technique}}}+\underbrace{A^{p+q}[X,Z]}_{\in\Im A^{p+q}}+\underbrace{[A^pX,A^qZ]-A^q[A^pX,Z]}_{=0\text{ as }X\in\ker A^p}.$$


{\em Regularity.} If $A$ is of class $C^r$, the distributions $\Im A^p$ and $\ker A^q$ are of class $C^{r}$ {\em i.e.\@} the foliations are of class $C^{r+1}$, so are the coordinates. If $A$ is holomorphic, everything is.\hspace{\fill}$\Box$\medskip

\noindent{\bf Proof of the theorem.} If $A$ is integrable, the integrability of $\ker A^p$ and the nullity of ${\cal N}_A$ are immediate. Let us prove the converse. The proof, when directly written in the general case $A^{n-1}\neq A^n=0$ with $n$ any integer, is a cumbersome induction. So we state it in cases $n=2$ and $n=3$, where all the arguments are involved. Then we give the structure of the induction for the general case. We also suppose that $A$ is of class $C^\infty$ and postpone the remarks about regularity when $A$ is analytic or of class $C^r$.\medskip

Proving that $A$ is integrable amounts to building a  field of basis $\beta$ on $\VV$ such that:\medskip

{\bf (i)} $\Mat_\beta(A)$ is constant,\medskip

{\bf (ii)} any two vector fields of $\beta$ commute (in other terms, the field $\beta$ is integrable).\medskip

\noindent{\em Proof for $n=2$.} Here $\II\subset\KK\subset \VV$; take coordinates $((x^1_\beta)_\beta,(x^2_\beta)_\beta,(x^3_\beta)_\beta)$ of $\VV$, adapted to this flag of foliations. Set $(Z_i)_{i=1}^{d_1}$ the coordinate vectors $\left(\frac{\partial}{\partial x^2_\ast}\right)$ 
 and $(Z_i)_{i=d_1+1}^{d_1+d_2}$ the $\left(\frac{\partial}{\partial x^3_\ast}\right)$, so that $Z_i\in\ker A$ for $i\leqslant d_1$. Then $((AZ_i)_{i>d_1},(Z_i)_{i=1}^{d_1+d_2})$ is a basis field on $\VV$, the $(Z_i)_{i=1}^{d_1}$ belonging to $\ker A$. Thus, in this basis:
$$\operatorname{Mat}(A)={\rm constant}=\left(\begin{array}{ccc}
0&0&I_{d_2}\\
0&0&0\\
0&0&0\end{array}\right).$$
We now replace the $Z_i$ by some {\em commuting} $Z'_i$, letting the form of $\Mat(A)$ unchanged. All vector fields are $\pi$-basic ({\em i.e.\@} $\II$-basic), so all brackets are in $\Im A$ (the fields ``commute modulo $\Im A$''). Moreover, as $\NN_A=0$, for any $i,j$, $[AZ_i,AZ_j]=A[Z_i,AZ_j]+A[AZ_i,Z_j]-A^2[Z_i,Z_j]\in\Im A^2=\{0\}$, so the $AZ_i$ commute. 
 Let $\SSS$ be the level $\{x^1=0\}$ (transverse to $\II$) and $\Phi^t_j$ be the flow of $AZ_j$, for $j>d_1$. Those flows commute and define a diffeomorphism $\Phi:(m,(t_j)_{j=d_1+1}^{d_1+d_2})\mapsto\Phi^{t_{d_1+1}}_{d_1+1}\circ\ldots\circ\Phi^{t_{d_1+d_2}}_{d_1+d_2}$ from $\SSS\times B_{\K^{d_2}}(\varepsilon)$ on a neighbourhood of $\SSS$ in $\VV$. We now set $(Z'_i)_i:=(Z_i)_i$ along $\SSS$, and {\em push} them by the flows $\Phi^t_j$. Formally: $Z'_i(\Phi(m,(t_j)_{i=d_1+1}^{d_1+d_2}))=\dd(\Phi^{t_{d_1+1}}_{d_1+1}\circ\ldots\circ\Phi^{t_{d_1+d_2}}_{d_1+d_2})(m).Z_{i}(m).$ Then: {\bf (a)} the $Z'_i$ are coordinate vector fields along $\SSS$, and are pushed forward on $\VV$ by commuting flows, so they commute everywhere, and by construction they commute with the fields $AZ_i$ (apply the Jacobi identity); {\bf (b)} the flows $\Phi^t_j$ respect the leaves of $\II$, so $Z_i-Z'_i\in\Im A$, so $AZ_i=AZ'_i$; {\bf (c)} the flows $\Phi^t_j$ respect the leaves of $\KK$. To check {\bf (c)}, take $Z$ a vector field in $\ker A$, then as $\NN_A=0$, $[AZ_j,Z]\in\ker A$: $A[AZ_j,Z]=[AZ_j,AZ]+A^2[Z_j,Z]-A[Z_j,AZ]=0$.

Let us conclude. By {\bf (a)} and {\bf (b)}, the basis field $\beta=((AZ'_i)_{i>d_1},(Z'_i)_{i=1}^{d_1+d_2})$ consists of commuting vector fields. By {\bf (c)}, the $(Z'_i)_{i=1}^{d_1}$, obtained by pushing the $(Z_i)_{i=1}^{d_1}$ by the $\Phi^t_j$, belong still everywhere to $\ker A$, so $\Mat_\beta(A)$ is unchanged. We are done.\medskip

\noindent{\em Proof for $n=3$.} We see here that in general, we will need an {\em induction}. This time, let $(x^{u,v}_\beta)$ be a coordinate system given by Lemma \ref{section} for the array of foliations we deal with:
$$\begin{array}{cccccc}
\II^2\\
\cap\\
(\II^1\cap\KK)&\subset&\II^1\\
\cap&&\cap\\
\KK&\subset&\KK^2&\subset&\VV.
\end{array}\quad\text{\begin{tabular}{l}(So the coordinates\\parametrising the\\ foliations are,\\ correspondingly:)\end{tabular}} 
\begin{array}{ccc}
(x^{1,1}_\beta)_\beta\\
(x^{2,1}_\beta)_\beta&(x^{2,2}_\beta)_\beta\\
(x^{3,1}_\beta)_\beta&(x^{3,2}_\beta)_\beta&(x^{3,3}_\beta)_\beta.
\end{array}$$
Let us set  $(Z_i)_{i=1}^{d_1+d_2+d_3}$ the coordinate vectors $\frac{\partial}{\partial x^{3,\star}_\star}$, in a way such that $Z_i\in\ker A$ for $i\leqslant d_1$ and $Z_i\in\ker A^2$ for $i\leqslant d_1+d_2$ {\bf (a)}. The fields $Z_i$ are $\KK^q$-basic for all $q$ {\bf (b)}; in other words, for any $j$, as soon as $Z_i\in\ker A^q$, $[Z_j,Z_i]\in\ker A^q$. By construction, the $Z_i$ are also $\pi$-basic, so for any $p$, $[Z_i,A^pZ_j]\in\Im A$ {\bf (c)}. The family $((A^2Z_i)_{i>d_1+d_2},(AZ_i)_{i>d_1},(Z_i)_{i})$ is a basis field on $\VV$, consisting of vector fields commuting modulo $\Im A$ and in which, because of {\bf (a)}:
$$\operatorname{Mat}(A)={\rm constant}=M=\left(\begin{array}{cccccc}
0&0&I_{d_3}&0&0&0\\
0&0&0&0&I_{d_2}&0\\
0&0&0&0&0&I_{d_3}\\
0&0&0&0&0&0\\
0&0&0&0&0&0\\
0&0&0&0&0&0\end{array}\right).$$
We now replace the $Z_i$ by some $Z'_i$ {\em commuting modulo $\Im A^2$}, letting the form of $\Mat(A)$ unchanged. As above, we take $\SSS$ the level $\{x^{2,\star}=0,x^{1,\star}=0\}$ (transverse to $\II$), and $(\Phi^t_\alpha)_{\alpha=1}^N$ the flows of the fields $(AZ_k)_{k>d_1},(A^2Z_k)_{k>d_1+d_2}$, arbitrarily indexed by some $\alpha\in\llbracket1,N\rrbracket$. Because of {\bf (c)} and as $\NN_A=0$, $[A^pZ_i,A^qZ_j]\in\Im A^2$ for $p\geqslant1$ and $q\geqslant 1$ {\bf (d)}, so those flows commute modulo $\Im A^2$. As in the case $n=2$, we build a diffeormorphism $\Phi:(m,(t_\alpha)_{\alpha=1}^{N})\mapsto\Phi^{t_{1}}_{1}\circ\ldots\circ\Phi^{t_{N}}_{N}$ from $\SSS\times B_{\K^{N}}(\varepsilon)$ on a neighbourhood of $\SSS$ in $\VV$; $\Phi$ depends on the arbitrary order of the $\Phi^t_\alpha$ but it does not matter. We define similarly the $Z'_i$. As the $\Phi^t_\alpha$ commute modulo $\Im A^2$, so do the $Z'_i$ with each other, and with the $AZ_j$ and $A^2Z_j$:
$$\text{for all $(i,j)$ and $p\geqslant1$, }\quad[Z'_i,Z'_j]\in\Im A^2\ \text{\bf (e)}\quad\text{and }[Z'_i,A^pZ_j]\in\Im A^2\ \text{\bf (f)}.$$
Besides, as the $\Phi^t_\alpha$ preserve each leaf of $\II$, for any $i$, $Z'_i\equiv Z_i [\Im A]$. Thus for $p\geqslant1$, $A^pZ'_i\equiv A^pZ_i\ [\Im A^2]$ {\bf (g)}. For $p\geqslant1$ and $q\geqslant 1$, it comes from {\bf (f)} and {\bf (g)}: $[Z'_i,A^pZ'_j]\in\im A^2$, and from {\bf (d)} and {\bf (g)}: $[A^pZ'_i,A^qZ'_j]\in\im A^2$. So the basis field $((A^2Z'_i)_{i>d_1+d_2},(AZ'_i)_{i>d_1},(Z'_i)_{i})$ is made of vector fields commuting modulo $\Im A^2$.

We also check that the $Z'_i$ still satisfy {\bf (a)} and {\bf (b)} {\em i.e.\@} that, if $Z_i\in\ker A^q$, then $Z'_i\in\ker A^q$, and that the $Z'_i$ are $\KK^q$-basic. The $Z'_i$ are equal to the $Z_i$ along $\SSS$, and are pushed by the flows $\Phi^t_\alpha$ of the $A^pZ_j$, $p\geqslant 1$. The wanted results follow from the fact that those flows preserve each foliation $\KK^q$ ---equivalently, the fact that the $A^pZ_j$ are $\KK^q$-basic. Indeed, take $Z$ any vector field in $\ker A^q$, then, as $\NN_{A^p,A^q}=0$, $[A^pZ_j,Z]\in\ker A^q$:
$$A^q[A^pZ_j,Z]=[A^pZ_j,\underbrace{A^qZ}_{\text{\tiny$=0$}}]+A^{p+q}\underbrace{[Z_j,Z]}_{\text{\tiny$\in\ker A^q$ by {\bf (b)}}}-A^p[Z_j,\underbrace{A^qZ}_{\text{\tiny$=0$}}]=0.$$

{\em Results at that step, and a remark:} we just get a basis field $((A^2Z'_i)_{i>d_1+d_2},\linebreak[2](AZ'_i)_{i>d_1},\linebreak[2](Z'_i)_{i})$, made of vector fields commuting modulo $\Im A^2$, and satisfying {\bf (a)} and {\bf (b)}. Consequently, $\Mat(A)$ in this basis has the constant writing $M$ given above. Using again $\NN_A=0$, we get moreover that, for $p\geqslant 1$ and $q\geqslant1$, $[A^pZ'_i,A^qZ'_j]\in\Im A^3=\{0\}$.\medskip

{\em End of the proof.} Using the remark just above, and iterating the process, we get new flows $\Phi^{\prime t}_\alpha$ which, this time, commute. It comes a new $\Phi'$, and new fields $Z''_i$ which, this time, commute with each other and with the $A^pZ'_i$, $p\geqslant 1$. We conclude as for the case $n=2$.\medskip

\noindent{\em Proof for any $n$.} The case $n=3$ contains all arguments. So here we only state the structure of the induction. We set $d_a=\dim(\dd\pi(\ker A^a)/\dd\pi(\ker A^{a-1}))$. {\em Remark:} with this notation, the invariant factors of $A$ are $(\underbrace{X,\ldots,X}_{\text{\tiny$d_1$ times}},\underbrace{X^2,\ldots,X^2}_{\text{\tiny$d_2$ times}},\ldots,\underbrace{X^{n},\ldots,X^{n}}_{\text{\tiny$d_n$ times}})$.

Take a coordinate system $(x^{u,v}_\beta)$ as given by Lemma \ref{section} for the array of foliations displayed on page \pageref{tableau_feuilletages}. 
Let us set  $(Z^{(0)}_i)_{i=1}^{d_1+\ldots+d_n}$ the coordinate vectors $\frac{\partial}{\partial x^{n,\star}_\star}$, in a way such that $Z_i\in\ker A$ for $i\leqslant d_1$, $Z_i\in\ker A^2$ for $i\leqslant d_1+d_2$ {\em etc.} We set $\SSS$ the level $\{x^{1,\star}=0,\ldots,x^{{n-1},\star}=0,\}$, transverse to $\II$.

 This builds vector fields $Z_{i}^{(0)}$ satisfying the following induction hypothesis, at step $k=0$:
$$\text{{\bf\mathversion{bold}(H$_k$)\mathversion{normal}} }\left\{\begin{array}{l}\text{\bf (1) } \text{the }Z_{i}^{(k)}\ \text{are equal to the $Z^{(0)}_i$ along $\SSS$,}\\
\text{\bf (2) } \text{the }Z_{i}^{(k)}\ \text{are $\KK^q$-basic for all $q$,}\\
\text{\bf (3) } \text{For any $q$ and for }i\leqslant d_1+\ldots+d_q,\ Z_{i}^{(k)}\in\ker A^q,\\
\text{\bf (4) } \text{for any $(a,b)\in\N^2$, for any $(i,j)$, }[A^aZ_{i}^{(k)},A^bZ_{j}^{(k)}]\in\Im A^{k+1},\\
\text{\bf (5) } \text{for any $(a,b)\in(\N^\ast)^2$, for any $(i,j)$, }[A^aZ_{i}^{(k)},A^bZ_{j}^{(k)}]\in\Im A^{k+2}.
\end{array}\right.$$

Setting $\beta^{(k)}=(A^pZ_i^{(k)})_{i,p}$, it follows from {\bf\mathversion{bold}(H$_k$)\mathversion{normal}} that $\beta^{(k)}$ is a basis field on $\VV$, in which $\Mat(A)$ has a constant Jordan form, of the type given for the case $n=3$.

If fields $Z_{i}^{(k)}$ are built, satisfying {\bf\mathversion{bold}(H$_k$)\mathversion{normal}}, then you introduce the flows $(\Phi^{(k)t}_\alpha)_{\alpha=1}^N$ of the fields $A^pZ^{(k)}$ for $p\geqslant1$, you set $Z_{i}^{(k+1)}=Z^{(0)}_i$ along $\SSS$ and then push the $Z_{i}^{(k+1)}$ on the whole $\VV$ by the flows $\Phi^{(k)t}_\alpha$, in an arbitrary order. The very arguments given for the case $n=3$ show that the $Z_{i}^{(k+1)}$ satisfy {\bf\mathversion{bold}(H$_{k+1}$)\mathversion{normal}}. The induction propagates.\medskip

\noindent {\em Conclusion and regularity questions.} If $A$ is of class $C^\infty$, the basis field $\beta^{(n-1)}$ consists of commuting vector fields, so we are done. If $\K=\R$ and $r=\omega$, or $\K=\C$, the flows $\Phi^t_i$ are given by the Cauchy-Kovalevskaya theorem, so all remains analytic and we are also done.

In case $A$ is only of class $C^r$, $r<\infty$, each step of the induction loses one order of regularity. Indeed, if $\Phi^t_i$ is the flow of some $A^pZ^{(k-1)}_i$ of class $C^r$, $\Phi^t_i$ is also $C^r$, so the $Z^{(k)}_j$, defined as the $Z^{(k-1)}_j$ pushed by the {\em differential} of the $\Phi^t_i$, are only $C^{r-1}$. So we may lose $n$ orders of regularity. Modifying slightly the end of the proof, we see that we lose only $n-1$.

Carrying on the induction up to {\bf\mathversion{bold}(H$_{n-1}$)\mathversion{normal}} would provide some $C^{n-r+1}$ fields $Z^{(n-1)}_{i}$, but as possibly $n-r+1=0$, this is useless: commuting fields of class $C^0$ but not $C^1$ do not provide corresponding coordinate functions, in general. Instead, we use directly the $C^{r-n+2}$-diffeomorphism $\Phi:(m,(t_i)_{i=1}^N)\mapsto \Phi^{t_N}_N\circ\ldots\circ\Phi^{t_1}_1(m)$ of this $(n-1)^{\text{th}}$ step of the induction. As the fields $A^aZ^{(n-2)}_{i}$, with $a>0$, commute, and parametrising $m\in\SSS$ by its coordinates $(x^{n,i}_{\alpha})_{i,\alpha}$, $\Phi$ is nothing but a local parametrisation of $\VV$ by a system of coordinates of class $C^{n-r+2}$, with coordinate vectors all the $A^aZ^{(n-1)}_{i}$ with $a\geqslant 0$. These coordinate vectors form a basis field in which $\Mat(A)$ has a constant Jordan form. We are done.

Eventually, the condition that  $A$ is of class $C^{n-1}$ is necessary, and the given  regularity of the integral coordinates is optimal: this follows from Example \ref{exemple_regularite} in the next section.\hfill$\Box$\\

\noindent{\bf Proof of Corollary \ref{corollaire}.} The integrability of the characteristic subspaces of $A$ amounts to their involutivity, through the Frobenius criterion. In turn this is implied by the nullity of ${\cal N}_A$. First, let us build integral coordinates on the integral leaf of each characteristic subspace, through the origin. On each characteristic subspace, take $A=S+N$ the ``semi-simple + nilpotent'' decomposition of $A$.

On the integral leaf of the spaces relative to some real eigenvalue $\lambda$, $S=\lambda\Id$, so applying Theorem \ref{theorem} to the nilpotent part $N$ gives the coordinates.

On the integral leaf of the other spaces, $S=\lambda\Id+\mu J$ for some $J$ with $J^2=-\Id$. For any commuting endomorphisms $U$ and $V$, $\NN_{U+V}=\NN_U+\NN_V+\NN'_{U,V}+\NN'_{V,U}$, so using Proposition \ref{technique}, we get that, for any $P\in\K[X]$, $\NN_{P(A)}=0$ as soon as $\NN_A=0$. So here $\NN_J=\NN_N=0$, $J$ is integrable by the integrability condition for complex structures, and $N$, viewed as a complex endomorphism, is integrable by Theorem \ref{theorem}.

Finally, take the unique ``product'' coordinate system extending the coordinates built above, on $\R^d$: it is integral for $A$. Indeed,  for each characteristic subspace $E$ of $A$, you may take $Q\in\R[X]$ such that $Q(A)_{|E}=A_{|E}$ on $E$ if $A_{|E}$ is invertible, $Q(A)_{|E}=A_{|E}+\Id_E$ on $E$ if $A_{|E}$ is nilpotent, and $Q(A)_{|F}=0$ on the sum $F$ of the other characteristic subspaces. To prove that the matrix of $A$ is constant in our coordinates, we must check that $(L_YA)(X)=0$, {\em i.e.} that $[Y,AX]=0$, for any coordinate vector fields $X$ tangent to $E$ and $Y$ tangent to $F$. Now, as $\NN'_{Q(A),A}=0$: $Q(A)[Y,AX]=[Q(A)Y,AX]-A[Q(A)Y,X]+Q(A).A[Y,X]=0$. As $[Y,AX]\in E$ and as $Q(A)_{|E}$ is invertible, we are done.\hfill$\Box$

\section{Some additional results and examples}\label{bonus}

\begin{prop}{\bf\mathversion{bold}[A higher partial regularity of the coordinates, when $\K=\R$]\mathversion{normal}}\label{regularite} In restriction to each integral leaf $\II^k$ of $\Im A^k$, for each $k\leqslant n-1$, the coordinates built by Theorem \ref{theorem} are of class $C^{r-n+2+k}$, and in general not more. Besides, The coordinates that are constant along the leaves of $\II^{k}$ are of class $C^{r+2-k}$, and in general not more. 
\end{prop}

\noindent{\bf Proof.} The optimality: ``in general not more'' follows from Example \ref{exemple_regularite} below. To prove the announced regularities, recall that the $(A^aZ^{(n-1)}_{i})_{a,i}$ are the coordinate vectors finally obtained in Theorem \ref{theorem}. In the proof of Theorem \ref{theorem}, each vector field $A^aZ^{(n-1)}_{i}$ is well-determined modulo $\Im A^{a+k+1}$ from the moment that $A^aZ^{(k)}_{i}$ is defined {\em i.e.\@} $A^aZ^{(k)}_{i}\equiv A^aZ^{(k+1)}_{i}\equiv \ldots \equiv A^aZ^{(n-1)}_{i}\ [\Im A^{a+k+1}]$. In particular:\medskip

{\bf (i)} The $(A^aZ^{(n-1)}_{i})_{a\geqslant n-k-1}$ are well-determined from step $k$ of the induction {\em i.e.\@} $A^aZ^{(n-1)}_{i}\linebreak[2]= A^aZ^{(k)}_{i}$. But the $Z^{(k)}_{i}$ are of class $C^{r-k}$, so the $(A^aZ^{(k)}_{i})_{a\geqslant n-k-1}$ are of class $C^{r-k}$.\medskip

{\bf (ii)} As the $Z^{(k)}_{i}$ are of class $C^{r-k}$, and as for any $a\geqslant 0$, $A^aZ^{(n-1)}_{i}\equiv A^aZ^{(k-1)}_{i}\ [\Im A^{k}]$ then the $[A^aZ^{(n-1)}_{i}\ \text{mod}\ \Im A^{n-k}]$ are all of class $C^{r+1-k}$.\medskip

\noindent Now, the $A^aZ^{(n-1)}_{i}$ with $a\geqslant n-k-1$ are the coordinate vectors along the leaves of $\II^{n-k-1}$. So by {\bf (i)}, the coordinates are of class $C^{r-k+1}$ along those leaves, the first claim. For the last claim,  denote the coordinates given by Theorem \ref{theorem} by $(y_i)_{i=1}^n=((y_{i,j})_{j=1}^{N_i})_{i=1}^n$, on such a way that the leaves of $\II^{k}$ are the levels of the $N$-tuple $(y_i)_{i>n-k}$. Then the $(\frac{\partial}{\partial y_i})_{i>n-k}$ are the $((A^aZ^{(n-1)}_{i})_{i})_{a<k}$ and we have to show that the $(y_i)_{i>n-k}$ are of class $C^{r+2-k}$ for any $k\geqslant1$. Take any coordinate system $(y'_i)_{i=1}^n=((y'_{i,j})_{j=1}^{N_i})_{i=1}^n$ {\em of class $C^{r+1}$} such that the leaves of $\II^{k}$ are the levels of the $N$-tuple $(y'_i)_{i>n-k}$. As the $(y_i)_{i=1}^n$ share the same property, the matrix $M=(\dd y_i(\frac{\partial}{\partial y'_j}))_{i,j=1}^n$ is upper block triangular, as well as $\Mat_{({\partial}/{\partial y'_j})_{j=1}^n}(\frac{\partial}{\partial y_i})_{i=1}^n=M^{-1}$. Thus, for each $k\geqslant1$, $\left(\dd y_i(\frac{\partial}{\partial y'_j})\right)_{i>n-k,j\leqslant n-k}=0$ and:
$$\left(\dd y_i\left(\frac{\partial}{\partial y'_j}\right)\right)_{i,j>n-k}=\left[\Mat_{({\partial}/{\partial y'_j})_{j>n-k}}\left(\frac{\partial}{\partial y_i}\ \text{mod}\ \Im A^{k}\right)_{i>n-k}\right]^{-1}.$$
By {\bf (ii)}, the matrix on the right side is of class $C^{r+1-k}$ so the $(\dd y_i)_{i>n-k}$ are of class $C^{r+1-k}$ and the  $(y_i)_{i>n-k}$ are of class $C^{r+2-k}$.\hfill$\Box$\medskip

For the next proposition, we introduce the following terminology.

\begin{vocabulary}\label{section_respectant_K} A section $\sigma$ of $\pi:\VV\rightarrow\VV/\II$ is said here to respect the foliations $\KK^1\subset\ldots\subset\KK^{n-1}$ if for all $p\in\N$, $\dd\sigma.(\dd\pi(\ker A^p))\subset\ker A^p.$ This amounts to saying that the image of $\sigma$ is the level $\{x^{1,\star}=0,\ldots,x^{{n-1},\star}=0,\}$ of some coordinate system given by Lemma \ref{section}.
\end{vocabulary}

\begin{prop}{\bf[Uniqueness of the integral coordinates]}\label{unicite} Let $A$ be an integrable field of nilpotent endomorphisms, $\II$ be the integral foliation of the distribution $\Im A$, $\KK^p$ be that of each $\ker A^p$ for $p\in\N$, and $\pi$ be the projection $\VV\rightarrow\VV/\II$. Then a system of integral coordinates for $A$, in which $\Mat(A)$ is a constant Jordan matrix, is uniquely given by the independent choice of:\medskip

-- a section $\sigma$ of $\pi$, respecting the foliations $\KK^1\subset\ldots\subset\KK^{n-1}$, in the sense of Vocabulary \ref{section_respectant_K},\medskip

-- coordinates $((\zbar_{1,\alpha})_{\alpha=1}^{d_1},\ldots,(\zbar_{n,\alpha})_{\alpha=1}^{d_n}))$ of $\pi(\VV)$ respecting the foliations $\pi(\KK^1)\subset\ldots\subset\pi(\KK^{n-1})$ {\em i.e.\@} such that the leaves of each $\pi(\KK^p)$ are the levels of  $\{((\zbar_{p+1,\alpha})_{\alpha},\ldots,(\zbar_{n,\alpha})_{\alpha}\}$.\medskip

\noindent More precisely, there is a unique ``Jordan'' coordinate system $((z'),(z_1),\ldots,(z_n))$ for $A$, characterised by the fact that:\medskip

-- $((z_1),\ldots,(z_n))=\pi^\ast(\zbar_1,\ldots,\zbar_n)$, (the levels of this $k$-tuple are the leaves of $\II$),\medskip

-- the coordinates $(z')$ are determined by the fact that $\{(z')=0\}$ is the image of $\sigma$ and that the $k$-tuple $(\frac{\partial}{\partial z'_i})_i$ is equal to that of the non null $(A^a\frac{\partial}{\partial z_{n,j}})_{n,j}$. (The coordinates $(z')$ parametrise the leaves of $\II$.) Explicitly, the fields of coordinate vectors are the $k$-tuple:
\begin{gather*}\textstyle\left((A^{n-1}\frac{\partial}{\partial z_{n,j}})_{j=1}^{d_n},((A^{n-2}\frac{\partial}{\partial z_{n-i,j}})_{j=1}^{d_{n-i}})_{i=1}^2,((A^{n-3}\frac{\partial}{\partial z_{n-i,j}})_{j=1}^{d_{n-i}})_{i=1}^3,\ldots,((\frac{\partial}{\partial z_{n-i,j}})_{j=1}^{d_{n}})_{i=1}^{n-i}\right).
\end{gather*}
\end{prop}

\noindent{\bf Proof.} We show that, once the image $\SSS$ of $\sigma$ and the fields $Z_{i}$ along it are chosen, the extension of the $Z_{i}$ satisfying Theorem \ref{theorem} is unique. Take $\widetilde{Z}_{i}$ another such extension. As $\widetilde{Z}_{i}=Z_{i}$ along $\SSS$ and as both fields are $\II$-basic, $\widetilde{Z}_{i}\equiv Z_{i}\ [\Im A]$ on $\VV$. So for $a>0$, $A^a\widetilde{Z}_{i}\equiv A^aZ_{i}\ [\Im A^{2}]$ on $\VV$. Now $[\widetilde{Z}_i,A^a\widetilde{Z}_i]=0$, so $[\widetilde{Z}_i,A^aZ_i]\in\Im A^2$. As $\widetilde{Z}_{i}=Z_i$ along $\SSS$, and as the saturation of $\SSS$ by the flows of the $(A^aZ_i)_{i,a>0}$ is the whole $\VV$, it comes that $\widetilde{Z}_{i}\equiv Z_i\ [\Im A^2]$ everywhere. By induction, $\widetilde{Z}_{i}\equiv Z_i\ [\Im A^k]$ for all $k$ and we are done.\hfill$\Box$

\begin{rem}\label{non_elliptique}Let us consider the particular case of an endomorphism field $A$, constant in the natural coordinates of the compact manifold $T=\R^d/\Z^d$. It follows from Proposition \ref{unicite} that the space of (global) integral coordinates for $A$ is infinite dimensional. This shows that such coordinates are not the solution of an elliptic problem. Instead, they appear naturally as the solution of a system of O.D.E., with an ``initial condition'' arbitrarily chosen in some infinite dimensional function space. This holds as soon as the minimal polynomial of $A$ contains a factor $(X-a)$, $a\in\R$, may $A$ be invertible or not.
\end{rem}

The following example shows that, in Theorem \ref{theorem}, $A$ has to be of class $C^{n-1}$, and that the integral coordinates may be not more regular than claimed in it and in Proposition \ref{regularite}.

\begin{example}\label{exemple_regularite} Consider $\R^n$ with coordinates denoted by $(x_i)_{i=1}^n$, take $r\in\N^\ast$ and $(\alpha)=(\alpha_i)_{i=1}^{n-1}$ an $(n-1)$-tuple of functions in $C^r(\R^n,\R)$, with $\alpha_{n-1}>0$ . Set $A=A_{(\alpha)}$ defined by $A_{(\alpha)}(\frac{\partial}{\partial x_1})=0$, $A_{(\alpha)}(\frac{\partial}{\partial x_i})=\frac{\partial}{\partial x_{i-1}}$ for all $i\in\llbracket 2,n-1\rrbracket$ and $A_{(\alpha)}(\frac{\partial}{\partial x_n})=\sum_{i=1}^{n-1}\alpha_i\frac{\partial}{\partial x_{i}}$. By construction, each $\ker A^p=\cap_{i=p+1}^{n}\dd x_i$ is integrable, and we check that $\NN_A=0$ if and only if:
$$\forall i,k\in\llbracket 1,n-2\rrbracket,\ \frac{\partial \alpha_i}{\partial x_k}=\frac{\partial \alpha_{i+1}}{\partial x_{k+1}}\ \text{and}\  \frac{\partial \alpha_i}{\partial x_1}=0.$$
We assume that this condition is satisfied, so Theorem \ref{theorem} applies. Notice that then, the knowledge of $\alpha_{n-1}$ determines all the other $\alpha_i$, up to a additive constant. Let us build the integral coordinates $(y_i)_{i=1}^n$ determined by an arbitrary choice of $\sigma$ and by the choice ``$\zbar_n=x_n$'' {\em i.e.\@} by $y_n=x_n$ (see Proposition \ref{unicite}). Notice that necessarily, $y_n=y_n(x_n)$, as the levels of both $x_n$ and $y_n$ are the integral leaves of $\ker A^{n-1}$. So, a reparametrisation $y_n(x_n)$ of the last coordinate amounts to multiply all the $\alpha_i$ by $1/y'_n(x_n)$, thus if $\alpha_{n-1}$ cannot by made independent of $x_n$ by a multiplication by some function of $x_n$ ---~we now assume this~---, we do not lose any generality by taking directly $y_n=x_n$. Similarly, a different choice of $\sigma$ amounts to add to each $x_i$ with $i\leqslant n-1$, some function $f(x_n)$. This lets all the $(\frac{\partial}{\partial x_i})_{i=1}^{n-1}$ unchanged and adds some linear combination of them to $\frac{\partial}{\partial x_n}$. In turn this lets the $\alpha_i$ unchanged, up to additive constants. So we do not lose generality.

Now the coordinates $y_i$ are determined by the above initial condition and the system:
$$M.A_{(\alpha)}.M^{-1}=A_{(0,\ldots,0,1)}\qquad \text{with }M=\left(\frac{\partial y_i}{\partial x_j}\right)_{i,j=1}^n.$$
As the $y_i$ must respect the foliations $\II^p$, notice that $y_k=y_k(x_k,\ldots,x_n)$. We let the reader check that the system, with the initial conditions, is equivalent to:
$$(\ast)\ \left\{\begin{array}[m]{l}
y_n=x_n\text{ and for }i\leqslant n-1,\ y_i=0\ \text{along }\{x_1=\ldots=x_{n-1}=0\},\\
\frac{\partial y_{n-k}}{\partial x_{n-1}}=\sum_{i=1}^kP_i\frac{\partial y_{n-k+i}}{\partial x_{n}}\quad\text{for all }k\in\llbracket 1,n-1\rrbracket,\\
\frac{\partial y_{k}}{\partial x_{n-i}}=\frac{\partial y_{k+i-1}}{\partial x_{n-1}}\quad\text{for all }k\in\llbracket 1,n-2\rrbracket\text{ and }i\in\llbracket 2,n-1\rrbracket,
\end{array}\right.$$
where the $P_i$ are the rational fractions in the $\alpha_i$ inductively defined by: $P_1=\frac1{\alpha_{n-1}}$ and $P_{i}=-\sum_{j=1}^{i-1}\frac{\alpha_{n-i+j-1}}{\alpha_{n-1}}P_j$. This system is overdetermined but, by Theorem \ref{theorem}, and as we have assumed that $\NN_A=0$, we now it is holonomic {\em i.e.\@} it admits a (here unique) solution. This solution is determined by the relation $y_{n-1}=P_1=\frac1{\alpha_{n-1}}$ and, by induction on $k$, by the equations, directly given by integration of $(\ast)$:
\begin{align*}
(\ast\ast)\quad y_{n-k}=&\int_0^{x_{n-1}}\sum_{i=1}^k(P_i\frac{\partial y_{n-k+i}}{\partial x_{n}})(x_1,\ldots,x_{n-2},t,x_n)\dd t\\+&\sum_{i=2}^{n-1}\int_0^{x_{n-i}}\frac{\partial y_{n-k+i-1}}{\partial x_{n-1}}(x_1,\ldots,x_{n-i-1},t,0,\ldots,0,x_n)\dd t.
\end{align*}
We have announced an effective example, so let us provide a simple one. Take $\alpha_{n-1}=1/(1+x_{n-1}\theta(x_n))$ with $\theta(t)=t^{r+1}$ for $t\geqslant 0$ and else $\theta(t)=0$. This $\alpha_{n-1}$ is of class $C^r$ and not of class $C^{r+1}$. This gives: $y_{n-1}=1+x_{n-1}\theta(x_n)$ and, by induction left to the reader:\medskip

-- for $k\leqslant r+1$, $y_{n-k}=\left(1+\frac{x_{n-1}^2}{2}\right)\theta^{(k-1)}(x_n)+z_{n-k}$, with $z_{n-k}$ of class $C^{r-k+2}$,\medskip

-- for $k\geqslant r+2$, $y_{n-k}$ is not defined.\medskip

 \noindent In Theorem \ref{theorem}, we want the $y_i$ to be of class (at least) $C^1$ ---~else writing $A$ in them makes no sense~---, so we must require here that $y_1$ is well defined and of class $C^1$ {\em i.e.\@} that $r\geqslant n-1$. Moreover, $y_1$ is of class $C^{r-n+2}$ and not of class $C^{r-n+3}$, so the regularity given in Theorem \ref{theorem} is optimal. Similarly, the example shows also the optimality of Proposition \ref{regularite}.
\end{example}

\begin{rem} We may add that if a vector field $V$ is of class $C^s$, its flow $\Phi^t_V$ is of class $C^s$ and, for a generic $V$, is not of class $C^{s+1}$. Thus if  $W$ is another vector field, of class $C^{s'}$ with $s'\geqslant s$, its image $(\Phi^t_V)_\ast W$ for $t\neq0$ is of class $C^{s-1}$ and, for a generic $V$, is not of class $C^{s}$. Used inductively in the proof of Theorem \ref{theorem}, this remark shows that, for a generic field $A$, the vector fields $Z^{(k)}_{i,j}$ are of class $C^{r-k}$ and {\em not more}. So for a generic $A$, the coordinates are not more regular than announced in Theorem \ref{theorem} and Proposition \ref{regularite}.
\end{rem}

The two little counter-examples \ref{noyau_non_integrable} and \ref{noyau_integrable} ensure the independance of both last conditions of Theorem \ref{theorem}.

\begin{example}\label{noyau_non_integrable}
Here is a field $A$ such that ${\cal N}_A=0$ and $\ker A$ is non involutive, with minimal nilpotence index of $A$ (2) and ambient dimension (4). In $\K^4$ with coordinates $(x_i)_{i=1}^4$, define $A$ by $A(\frac{\partial}{\partial x_1})=A(\frac{\partial}{\partial x_2})=0$, $A(\frac{\partial}{\partial x_3})=\exp(x_2)\frac{\partial}{\partial x_1}$ and $A(\frac{\partial}{\partial x_4})=\frac{\partial}{\partial x_1}$. All $[A^a\frac{\partial}{\partial x_i},A^b\frac{\partial}{\partial x_j}]$ for $\{a,b\}\subset\{0,1\}$ vanish except $[A\frac{\partial}{\partial x_3},\frac{\partial}{\partial x_2}]=-\exp(x_2)\frac{\partial}{\partial x_1}$, hence $\NN_A=0$. But $\ker A=\ker\alpha$ with $\alpha=\dd x_4+x_2\dd x_3$, and $\alpha\wedge\dd\alpha=\dd x_2\wedge\dd x_3\wedge\dd x_4\neq 0$ so $\ker A$ is not involutive.
\end{example}

\begin{example}\label{noyau_integrable}
Here is a field $A$ such that ${\cal N}_A\neq 0$ and $\ker A$ is involutive, with minimal ambient nilpotence index of $A$ (again 2) and dimension for it (again 4). Similarly, define $A$ by $A(\frac{\partial}{\partial x_1})=A(\frac{\partial}{\partial x_2})=0$, $A(\frac{\partial}{\partial x_3})=\exp(x_2)\frac{\partial}{\partial x_1}$ and $A(\frac{\partial}{\partial x_4})=\frac{\partial}{\partial x_2}$. All $[A^a\frac{\partial}{\partial x_i},A^b\frac{\partial}{\partial x_j}]$ for $\{a,b\}\subset\{0,1\}$ vanish except $[A\frac{\partial}{\partial x_3},\frac{\partial}{\partial x_2}]=[A\frac{\partial}{\partial x_3},A\frac{\partial}{\partial x_4}]=-\exp(x_2)\frac{\partial}{\partial x_1}$. So $\NN_A\neq 0$ as $\NN_A(\frac{\partial}{\partial x_3},\frac{\partial}{\partial x_4})=-\exp(x_2)\frac{\partial}{\partial x_1}$. But $\ker A=\ker(\dd x_3)\cap\ker(\dd x_4)$ is involutive.
\end{example}

\begin{rem} However, in Theorem \ref{theorem}, for some similarity types of endomorphisms $A$, the second condition may be omitted or relaxed, as it is (partially) implied by the first one. For instance, if $A$ is cyclic, then for every $p$, $\ker A^p=\im A^{n-p}$ is involutive. More generally, if for some $p$, $\dim(\ker A^p/\im A^{n-p})=1$, then $\ker A^p$ is involutive. Indeed, take $(Y_i)_i$ a basis field of $\Im A^{n-p}$ and $X$ a field such that $(X,(Y_i)_i)$ spans $\ker A^p$. As $\NN_A=0$, $[Y_i,Y_j]\in\Im A^{n-p}\subset \ker A^p$, besides $[X,X]=0$. Take $Z_i$ such that $Y_i=A^{n-p}Z_i$, then $A^p[X,Y_i]=-\NN'_{A^p,A^{n-p}}(X,Z_i)+[A^pX,A^{n-p}Z_i]-A^{n-p}[A^pX,Z_i]+A^n[X,Z_i]=0$ so $[X,Y_i]\in\ker A^p$, thus $\ker A^p$ is involutive.
\end{rem}

\begin{rem} If $A$ is nilpotent, ${\cal N}_A=0$ does not imply that the $\ker A^p$ are involutive. It gives however a weaker fact: if $X,Y\in\ker A^p$, then $[X,Y]\in\ker A^{2p}$. Indeed, by Proposition \ref{technique}, ${\cal N}_{A^p}(X,Y)=0$, so $A^{2p}[X,Y]=-[A^pX,A^pY]+A^p[X,A^pY]+A^p[A^pX,Y]=0$.
\end{rem}

\begin{rem}In Theorem \ref{theorem}, if $A$ is defined on $\VV=\K^d$, the integral coordinates may be in fact built on the whole $\K^d$. Indeed, Theorem \ref{theorem} builds coordinates on some ball $B(p,R_p)$, around any point $p$ of $\VV$, with $R_p$ depending only on the coefficients of the matrix $A$ around $p$, through the flows $\Phi_i$ appearing in the proof of the theorem. So on any precompact set of $\VV$, this $R_p$ is bounded from below by a positive constant. Now take any domain of the type $]-\alpha,\alpha[^d$, on which integral coordinates are defined; by what precedes and by the unicity result \ref{unicite}, these coordinates may be extended on some $]-\alpha',\alpha'[^d$ with $\alpha'>\alpha$. We are done.
\end{rem}

\begin{example-rem}A consequence of Corollary \ref{corollaire} is that, if $(M,\nabla)$ is a manifold with a torsion free affine connection, any parallel endomorphism field $A$ on $M$ is integrable. 

More generaly, an endomorphism field is integrable if and only if it is parallel for some torsion free affine connection $\nabla$ (compare \cite{libermann} Th.\@ 6.1). Indeed, if $A$ is integrable, define $\nabla$ by $\nabla\frac{\partial}{\partial v_i}=0$ in some integral coordinate system $(v_i)_{i=1}^d$. It is torsion free and immediately $\nabla A=0$. Conversely, suppose $\nabla A=0$ with $\nabla_UV-\nabla_VU=[U,V]$ for all vector fields $U$ and $V$. Then $\NN_A=0$ and, by the Frobenius criterion, each distribution $\ker A^p$ is integrable. Besides, $\nabla A=0$ implies that $A$ has constant invariant factors so Corollary \ref{corollaire} applies.

I do not know other significant examples where endomorphism fields satisfying  naturally the assumptions of Corollary \ref{corollaire} appear. 
\end{example-rem}

\end{document}